\documentclass[reqno]{amsart}
\usepackage{amssymb}  
\usepackage{amsmath} 
\usepackage{amsthm}  
\usepackage{bm} 
\usepackage{stmaryrd} 

\advance\hoffset by -1in
\advance\voffset by -1in

\newtheorem{prop}{Proposition}
\newtheorem{prob}{Problem}
\newtheorem{conj}{Conjecture}

\title[A Moufang loop with exceptional properties]{A Moufang loop with exceptional\\ properties of associators}
\author{Ilya B. Gorshkov, Alexandre N. Grichkov, Andrei V. Zavarnitsine}

\thanks{Supported by FAPESP (proc. 2014/13730-0), by the Russian Foundation for Basic Research
(project 13--01--00505)}
\date{}

\begin{document}
\begin{abstract} We construct a Moufang loop $M$ of order $3^{19}$ and a pair $a,b$ of its elements such that the set of all elements of $M$ that associate with $a$ and $b$ does not form a subloop. This is also an example of a nonassociative Moufang loop with a generating set whose every three elements associate.

{\sc Keywords:} Moufang loop, associator, subloop

{\sc MSC2010:}  20N05   

\end{abstract}
\maketitle

\section{Introduction}

A loop in which the identity $(xy)(zx)=(x(yz))x$ holds is called a Moufang loop. Basic properties of Moufang loops can be found in  \cite{pf}. For elements $x,y,z$ of a Moufang loop, we denote by $[x,y]$
the unique element $c$ such that $xy=(yx)c$ and by $(x,y,z)$ the unique element $a$ such that $(xy)z=(x(yz))a$.

The following two related problem arise naturally in the study of Moufang loops:

\begin{prob}\label{um} Let $L$ be a Moufang loop and let $a,b\in L$. Consider the set
$$
l_{a,b}=\{x\in L \mid (x,a,b)=1 \}.
$$
Is $l_{a,b}$ a subloop of $L$?
\end{prob}

\begin{prob}\label{dois} Let $L$ be a Moufang loop generated by a set $\{a_i\mid i\in I\}$. Suppose that
$(a_i,a_j,a_k)=1$ for all $i,j,k\in I$. Does it follow that $L$ is associative?
\end{prob}

For all known Moufang loops that we checked, both problems were solved in the affirmative. However,
we found a new example of a $4$-generator loop of order $3^{19}$ which refutes both assertions. This loop is
constructed in the next section. In relation to Problem \ref{dois}, we mention that the existence of 3-torsion in this counterexample is essential. We pose following conjecture.

\begin{conj}\label{conj} Then question in Problem \ref{dois} is answered in the affirmative if $L$ is a finite Moufang $p$-loop with $p\ne 3$.
\end{conj}

An evidence in favor of this conjecture is provided by the following simple fact from the theory of Malcev algebras.

\begin{prop}\label{prop} Let $M$ be a Malcev algebra over a field of characteristic other than $2,3$ with a generating set $X$. If $\ \mathrm{J}(x,y,z)=0$ for all $x,y,z\in X$ then $M$ is a Lie algebra.
\end{prop}
Indeed, every counterexample to this proposition would yield a counterexample to Conjecture \ref{conj}
which could be constructed using the Campbell--Hausdorff formula similarly to \cite{kuz}.

\section{The loop}

Let $M$ be a $19$-dimensional vector space over the field $\mathbb{F}_3$ of three elements.
Elements of $M$ will be written as tuples $x=(x_1,x_2,\ldots,x_{19})\in M$.
We introduce a new operation '$\circ$' on $M$ which is given,
for $x,y\in M$, by the formula
\begin{equation}\label{mult}
x\circ y = x + y + f,
\end{equation}
where $f=(f_1,\ldots,f_{19})$ and $f_k$ are polynomials in $x_i$, $y_j$ explicitly given below.

\begin{align*}
f_1   &= f_2  = f_3  = f_4 =  0,\\
  f_5   &=-x_2y_1,\
  f_6    =-x_3y_1,\
  f_7    =-x_4y_1,\
  f_8    =-x_3y_2,\
  f_9    =-x_4y_2,\
  f_{10} =-x_4y_3,\\
  f_{11}&=-x_2x_3y_1-x_2y_1y_3+x_5y_3-x_8y_1,\
  f_{12} =-x_2x_4y_1-x_2y_1y_4+x_5y_4-x_9y_1,\\
  f_{13}&=-x_3y_1y_2+x_6y_2+x_8y_1,\
  f_{14} =-x_3x_4y_1-x_3y_1y_4+x_6y_4-x_{10}y_1,\\
  f_{15}&=-x_4y_1y_2+x_7y_2+x_9y_1,\
  f_{16} =-x_4y_1y_3+x_7y_3+x_{10}y_1,\\
  f_{17}&=-x_3x_4y_2-x_3y_2y_4+x_8y_4-x_{10}y_2,\
  f_{18} =-x_4y_2y_3+x_9y_3+x_{10}y_2,\\  f_{19}&=-x_1x_2x_4y_3+x_1x_2y_3y_4+x_1x_3y_2y_4+x_1x_4y_2y_3-x_1y_2y_3y_4\\
        & -x_2x_3x_4y_1+x_2x_3y_1y_4+x_2x_4y_1y_3+x_3x_4y_1y_2-x_3y_1y_2y_4+x_1x_8y_4\\
        & -x_1x_9y_3+x_1x_{10}y_2-x_1y_2y_{10}+x_1y_3y_9-x_1y_4y_8-x_2x_6y_4+x_2x_7y_3-x_2x_{10}y_1\\
        & +x_2y_1y_{10}-x_2y_3y_7+x_2y_4y_6+x_3x_5y_4-x_3x_7y_2+x_3x_9y_1-x_3y_1y_9+x_3y_2y_7\\
        & -x_3y_4y_5-x_4x_5y_3+x_4x_6y_2-x_4x_8y_1+x_4y_1y_8-x_4y_2y_6+x_4y_3y_5.
\end{align*}

It can be checked that $(M,\circ)$ is a Moufang loop. The identity is the zero vector of $M$ and, for $x\in (M,\circ)$, we have
$$
x^{-1} = -x + h,
$$
where $h=(h_1,\ldots,h_{19})$ and the polynomials $h_k$ are as follows:
\begin{align*}
& h_1   =   h_2   =  h_3   =   h_4   =  0,\\
 & h_5   = -x_1x_2,\
   h_6   = -x_1x_3,\
   h_7   = -x_1x_4,\
   h_8   = -x_2x_3,\
   h_9   = -x_2x_4,\
   h_{10}= -x_3x_4,\\
 & h_{11}= -x_1x_8+x_3x_5,\
   h_{12}= -x_1x_9+x_4x_5,\
   h_{13}= x_1x_2x_3+x_1x_8+x_2x_6,\\
 & h_{14}= -x_1x_{10}+x_4x_6,\
   h_{15}= x_1x_2x_4+x_1x_9+x_2x_7,\
   h_{16}= x_1x_3x_4+x_1x_{10}+x_3x_7,\\
 & h_{17}= -x_2x_{10}+x_4x_8,\
   h_{18}= x_2x_3x_4+x_2x_{10}+x_3x_9,\
   h_{19}= -x_1x_2x_3x_4
\end{align*}

Let $e_1,\ldots,e_{19}$ be the standard basis of the original vector space $M$, i.\,e.,
$e_i=(\ldots,0,1,0,\ldots)$ with '$1$' at the $i$th place. Define
$$
a=e_1,\quad b=e_2,\quad c=e_3,\quad d=e_4.
$$
Then using the multiplication formula (\ref{mult}) we can check the following equalities in $(M,\circ)$:
\begin{align*}
& e_5=[a,b],\ e_6=[a,c],\ e_7=[a,d],\ e_8=[b,c],\ e_9=[b,d],\ e_{10}=[c,d],\\
& e_{11}=[[a,b],c],\ e_{12}=[[a,b],d],\ e_{13}=[[a,c],b],\ e_{14}=[[a,c],d],\\
& e_{15}=[[a,d],b],\ e_{16}=[[a,d],c],\ e_{17}=[[b,c],d],\ e_{18}=[[b,d],c],\\
& e_{19}=([a,b],c,d).
\end{align*}
Moreover, for any $n_1,\ldots,n_{19}\in \mathbb{Z}$, we have
$$
(\,\ldots((e_1^{n_1}\circ e_2^{n_2}) \circ e_3^{n_3})\ \ldots\ \circ e_{19}^{n_{19}}) = (\overline{n}_1,\ldots,\overline{n}_{19}),
$$
where $[n\mapsto \overline{n}]$ denotes the natural epimorphism $\mathbb{Z}\to \mathbb{F}_3$. In particular, $(M,\circ)=\langle a,b,c,d \rangle$. In this loop, the following equalities hold:
$$(a,b,c)=(a,b,d)=(a,c,d)=(b,c,d)=1, \qquad ([a,b],c,d)=e_{19}\ne 1$$
In particular, the loop is not associative, which answers in the negative Problem \ref{dois}. Moreover,
we have $a,b\in l_{c,d}$ and $[a,b]\not\in l_{c,d}$, which implies that $l_{c,d}$ is not a subloop and answers in the negative Problem \ref{um}.

For a definition and basic properties of Malcev algebras, see, e.\,g., \cite{sag}. We recall that for elements $a,b,c$ in a Malcev algebra one defines their Jacobian $\mathrm{J}(a,b,c)=(ab)c+(bc)a+(ca)b$.
We now prove Proposition \ref{prop}.

\begin{proof} Let $\mathrm{J}(M)$ be the ideal of $M$ generated by $\mathrm{J}(a,b,c)$ for all $a,b,c\in M$.
It suffices to show that $\mathrm{J}(M)=0$. Since $\mathrm{J}(a,b,c)$ is linear in each argument,
without loss of generality we may assume that $a,b,c$ are (nonassociative) words in $X$. We proceed by
induction on $n=|a|+|b|+|c|$, where $|w|$ denotes the length of a word $w\in M$. If $n=3$ then $\mathrm{J}(a,b,c)=0$ by assumption. If $n>3$ then we may assume that $|a|>1$ and $a=a_1a_2$ with
$|a|=|a_1|+|a_2|$. Using the identity \cite[(2.15)]{sag}
$$3\mathrm{J}(wx,y,z) = \mathrm{J}(x,y,z)w - \mathrm{J}(y,z,w)x - 2\mathrm{J}(z,w,x)y + 2\mathrm{J}(w,x,y)z$$
which holds in Malcev algebras over fields of characteristic
distinct from $2,3$, we have $\mathrm{J}(a,b,c)=0$ by induction.
\end{proof}

\end{document}